\documentclass[12pt]{article}

\usepackage{epsfig}
\usepackage{fullpage}
\usepackage{color}
\usepackage{amsmath}
\usepackage{setspace}
\usepackage{amsfonts}
\usepackage{amssymb}
\usepackage{graphicx}
\usepackage{hyperref}

\newtheorem{theorem}{Theorem}[section]
\newtheorem{lemma}[theorem]{Lemma}

\begin{document}

\title{Radial Symmetry of Large Solutions of Semilinear Elliptic Equations with Convection}
\author{ Ehsan Kamalinejad \thanks{Partially supported by NSERC 311685} \quad  and \quad Amir 
Moradifam \thanks{Supported by a MITACS Strategic Postdoctoral Fellowship
 }
\\
\small Department of Mathematics,
\small University of Toronto, \\
\small {\tt ehsan@math.toronto.edu} \\
\small {\tt amir@math.toronto.edu}
\\
}
\maketitle

\begin{abstract} 
We study radial symmetry of large solutions of the semi-linear elliptic problem $ \Delta u + \nabla h \cdot \nabla u = f(|x|,u)$, and we provide sharp conditions under which the problem has a radial solution. The result is independent of the rate of growth of the solution at infinity.\\


\end{abstract}


\section{Introduction}

Radial symmetry of the solutions of $\Delta u = f(|x|,u)$ on $\mathbb{R}^n$ is a well-studied problem and various conditions on the rate of growth and monotonicity of $f(|x|,u)$ as well as behaviour of $u(x)$ at infinity have been presented to guarantee radial symmetry of the solutions. In this paper we study radial symmetry of large solutions of the semi-linear elliptic problem
\begin{align} \label{equ-main}
\begin{cases}
\Delta u(x) + \nabla h(x) \cdot \nabla u(x)= f(|x|,u(x)) ~~~~~~&  x \in \mathbb{R}^n ~~~ (n\geq 2),\\
u(x)\longrightarrow \infty  & x\rightarrow \infty .
\end{cases}
\end{align}
We assume that for large values of $|x|$ and $u$ the function $f(|x|,u)$ is positive and superlinear, and that $\lim_{|x| \to \infty} u(x)= \infty$ but we do not assume a particular rate of growth at infinity for the solution. Our main focus is the effect of the convection term in radial symmetry of the solutions. The case $h \equiv 0$ with similar setting has been studied in \cite{T1} and \cite{T2}, and in contrast to the large boundary condition, symmetry of the small solutions $\lim_{|x| \to \infty}u(x)=0$ of the same problem has been studied in \cite{LY}, \cite{L1}, \cite{LN1} and \cite{N1}.\\   

If $u(x)$ is radial then all of the terms in (\ref{equ-main}) except perhaps $h(x)$ will be radial which automatically implies radial symmetry of $h(x)$ at least whenever $u$ is not constant. Thus, it is natural to assume that $h$ is radial and whenever clear, we abuse the notation $h(x)=h(|x|)$. We ask for the convection term $h$ to satisfy a particular integrability condition given by 
\begin{equation}\label{growth}
\exists R \geqslant 0 ~~,~~ \int_{R}^{\infty} e^{-h(r)} r^{1-n} dr <\infty.
\end{equation}
This condition is shown to be sharp in the sense that if violated, while all the other conditions hold, there are examples with no radial solution. Having this condition on $h$, a change of variable is proved to be well-defined which converts the radial solutions of the PDE into the solutions of a corresponding ODE. Then the available ODE theory developed in \cite{T1} combined with comparison arguments can be used to prove existence and symmetry of the solutions.\\


\section{Statements and Proofs}

To set the appropriate conditions on $f(|x|,u)$ we compare it with a function $g(r,s)$ that satisfies the following conditions:
\begin{itemize}
\item[($c1$)] $g(r,s)$ and $g_s(r,s)$ are continuous and positive on $\Omega=\lbrace (r,s) ~|~ r>r_0,s>s_0 \rbrace$ where $r_0$, $s_0$ are positive constants.
\item[($c2$)] $g(r,s)$ is superlinear in $s$ on $\Omega$ in the sense there exist $\lambda >1$ such that that $g(r,vs)\geqslant v^{\lambda}g(r,s)$ for all $v>1$ and $(r,s) \in \Omega$.
\item[($c3$)] $p(r)e^{h(r)}g(r,s)$ is monotone in $r$ on $\Omega$, where the function $p(r)$ is given by 
\begin{equation}\label{equ:p}
p(r):= -\int_{r}^{\infty} e^{-h(z)} z^{1-n} dz.
\end{equation}
\end{itemize}

The following theorem is the main result of this paper. 

\begin{theorem}\label{theo:radi}
Let $h(r)$ be continuous and satisfy (\ref{growth}). Let $f(r,s)$ and $f_s(r,s)$ be continuous and positive. Assume that there exist a function $g(r,s)$ such that
\[ \lim_{(r,s)\to(\infty, \infty)} \frac{f(r,s)}{g(r,s)}=1 \]
where $g(r,s)$ satisfies ($c1$-$c3$). Assume also that $f(r,s)$ is superlinear in $s$ on $\Omega$. Then 
\begin{itemize}
\item[(i)] All $C^2$ solutions of Problem (\ref{equ-main}) are radial.
\item[(ii)] If (\ref{equ-main}) has a $C^2$ solution then $\exists R \geqslant 0 ,~ \exists \hat{u}>0$ such that
\begin{equation}\label{abc}
-\int_{|x|>R} p(|x|)e^{h(|x|)}f(x,s)dx < \infty ~~~~ \forall s>\hat{u}
\end{equation}
where $p(r)$ is given by (\ref{equ:p}).
\item[(iii)] If in addition $f(|x|,u)$ satisfies ($c3$), Condition (\ref{abc}) is also a sufficient condition for existence of a solution to Problem (\ref{equ-main}).
\end{itemize}
\end{theorem}

In \cite{T1} and \cite{T2} Taliaferro studies relevance of the conditions on $f(r,s)$ for the problem without the convection term. For example it is shown that superlinearity of $f(r,s)$ is a sharp condition for radial symmetry of the solutions of \ref{theo:radi}. Indeed there are non-radial solutions of the problem \ref{theo:radi} when this condition fails.\\

Condition (\ref{growth}) on $h$ is a sharp condition in the sense that if it does not hold, then there are cases with no radial solution to Problem (\ref{equ-main}). To see this let $h(x)=\beta \log(|x|)$. Notice that Condition (\ref{growth}) holds for $\beta > 2-n$ and it is violated if $\beta \leqslant 2-n$. Consider the critical case when $\beta =2-n$ and let $f(r,s)=f(s)$ be a superlinear function. We claim that there is no radial solution to (\ref{equ-main}). Assume for the contrary that there exists a radial solution $u(x)=u(|x|)$. We have
\begin{align*}
f(u)&=\Delta u(x) + \nabla  log (|x|^{(2-n)})\cdot \nabla u(|x|)\\
&= \Delta u(x) + (2-n) \dfrac{x}{|x|^2} \cdot \frac{x}{|x|} u'(|x|)\\
\end{align*}
Hence
\begin{align*}
f(u) & = \left\lbrace u''(r) + \frac{n-1}{r}u'(r) \right\rbrace + \frac{2-n}{r}u'(r)\\
& = u''(r) + \frac{1}{r}u'(r)
\end{align*}
where $r=|x|$. Now define the radial function $v:\mathbb{R}^2 \longrightarrow \mathbb{R}$ by $v(y):=u(|y|)$. Then $v$ is a radial solution of the problem $\Delta v = f(v)$ in $\mathbb{R}^2$. This is a contradiction because Osserman showed in \cite{O} that for a superlinear function $f(v)$ the problem  $\Delta v(x) = f(v)$ has no large solution in $\mathbb{R}^2$. Therefore (\ref{equ-main}) has no radial solution or it has no solution at all, which both indicate necessity of the condition (\ref{growth}).\\


Based on the condition (\ref{growth}) on $h(x)$ we can use the following change of variables to transform radial solutions of (\ref{equ-main}) into the corresponding ODE solutions.
 
\begin{lemma}\label{theo:ode}
Let $h(r)$ satisfy (\ref{growth}) and let $f(|x|,u)$ satisfy the conditions of Theorem \ref{theo:radi}. Then $u(x)$ is a radial solution of Equation (\ref{equ-main}) if and only if $z(t):=u(p^{-1}(t))$ solves 
\begin{equation}\label{equ:odee}
  \begin{cases}
   z''(t)=F(t,z(t)), \\
   \lim_{t \to 0^-}z(t)=\infty,
  \end{cases}
\end{equation}
where $p(r)$ is given by (\ref{equ:p}), and $F(t,z)$ is given by 
\begin{equation}\label{equ:F}
F(t,z):=(p^{-1}(t))^{2n-2}e^{2h(p^{-1}(t))}f(p^{-1}(t),z).
\end{equation}
\end{lemma}

\textbf{Proof.} Let $r=|x|$, $t=p(r)$, and $z(t)=u(p^{-1}(t))$. This is a valid change of variable because by definition $p(r)$ is continuous and strictly increasing. We have
\begin{align*}
f(r,u(r))&=[u''(r)+\frac{n-1}{r}u'(r)]+h'(r)u'(r)\\ 
&=p'(r)^2 z''(p(r))+(p''(r)+\frac{n-1}{r}p'(r)+h'(r)p'(r))z'(p(r))\\ &=p'(r)^2 z''(p(r)).
\end{align*}
Therefore
\begin{align*}
z''(p(r))&=\dfrac{1}{p'(r)^2}f(r,z(p(r))\\
&=e^{2h(r)}r^{2n-2}f(r,z(p(r))\\
&=F(p(r),z(p(r))).
\end{align*}
Also the boundary condition $\lim_{|x| \to \infty }u(x)=\infty$ is equivalent to $\lim_{t \to 0^-}z(t)=\infty$ because $\lim_{r \to \infty} p^{-1}(r)=0$. \hfill $\square$\\

\textbf{Remark.} Note that the definition of $F(t,z)$ implies that for large values of $t$ and $z$ both $F(t,z)$ and $F_z(t,z)$ are continuous and non-negative, and that $F(t,z)$ is superlinear in $z$. This fact is useful when we study the ODE which corresponds to Equation (\ref{equ-main}). Lemma \ref{theo:ode} plays an important role in our arguments. In particular, in the proof of Theorem \ref{theo:radi} we need to construct two sequences of radial functions for the comparison arguments. The sequences can be constructed by the help of Lemma \ref{theo:ode} from the ODE counterparts described in Lemma \ref{lem:v} as follows. Assuming the conditions of Lemma \ref{theo:ode} hold, for each $M,m>s_0$ and $r_1>r_0$ there exists an increasing sequence $\lbrace \rho_k \rbrace^{\infty}_{k=1} \subseteq(r_1,\infty)$ with $\lim_{k \to \infty} \rho_k=\infty $, and two sequences of $C^2$ radial functions $\lbrace u_k(x)\rbrace$  and $\lbrace U_k(x) \rbrace$ that
\begin{itemize}
\item[(i)] $U_0(x)$ and $u_0(x)$, $u_1(x)$,... are radial solutions of  
\begin{align*}
  \begin{cases}
  \Delta u(x) + \nabla h(x) \cdot \nabla u(x) =f(|x|,u(x)) ~~~~~~& |x|\geqslant r_1, \\
   u(x)=m & |x|=r_1,
  \end{cases}
\end{align*}
\item[(ii)]$U_1(x)$, $U_2(x)$,... are solutions of  
\begin{align*}
  \begin{cases}
   \Delta U_k (x) + \nabla h(x) \cdot \nabla U_k (x)= f(|x|,U_k(x)) ~~~~~~& r_1 \leqslant |x|\leqslant \rho_k, \\
   U_k \longrightarrow \infty & |x| \longrightarrow \rho_k^- , \\
   U_k(x)=M & |x|=r_1,
  \end{cases}
\end{align*}
\item[(iii)]$\lim_{|x| \to \infty} u_0(x) =\lim_{|x| \to \infty} U_0(x)=\infty$,
\item[(iv)]$u_1(x)$, $u_2(x)$,... are all bounded as $|x| \to \infty$,
\item[(v)] For each $|x|>r_1$ we have $\lim_{k \to \infty}u_k(x)=u_0(x)$, and $\lim_{k \to \infty}U_k(x)=U_0(x)$.
\end{itemize}


We are now ready to prove Theorem \ref{theo:radi}.\\

\textbf{Proof of Theorem \ref{theo:radi}(iii).} We start by proving part (iii) where we have additional monotonicity condition ($c3$) on $f(|x|,u)$. In fact, we prove that, having (c3), Condition (\ref{abc}) is both necessary and sufficient for existence of a solution to (\ref{equ-main}). This fact will be useful in the proof of other parts. Let $t=p(|x|)$ where $p(r)$ is given by (\ref{equ:p}). We have
\begin{equation*}
-\int_ \Omega p(|x|)e^{h(|x|)}f(|x|,s)dx=
-\sigma_n \int_ R^\infty r^{n-1} p(r)e^{h(r)}f(r,s)dr,
\end{equation*}
where $\sigma_n$ is the perimeter of the unit ball in $\mathbb{R}^n$. Therefore
\begin{equation} \label{klk}
\begin{aligned}
- \int_\Omega p(|x|)e^{h(|x|)}f(x,s)dx&=
- \sigma_n \int_ R^\infty r^{2n-2} p(r)e^{2h(r)}f(r,s)(e^{-h(r)}r^{(1-n)})dr\\
&=-\sigma_n \int_ R^\infty p(r)F(p(r),s)p'(r) dr\\
&=-\sigma_n \int_{t_0}^0 tF(t,s)dt,
\end{aligned}
\end{equation}
where in the second equality we used the fact that $p'(r)=e^{-h(r)}r^{(1-n)}$. Assuming that (\ref{abc}) holds, (\ref{klk}) implies that
\begin{equation} \label{gk}
-\int_{t_0}^0 tF(t,s)dt<\infty ~~~\forall s>\hat{u}.
\end{equation}
By the ODE lemma \ref{lem:integ} the condition (\ref{gk}) is a necessary and sufficient condition for existence of a solution $z(t)$ of $z''(t)=F(t,z(t))$. By Lemma \ref{theo:ode} the solution $z(t)$ of the ODE $z''(t)=F(t,z(t))$ can be transformed into a radial solution $u(x)=z(p(|x|))$ of the equation (\ref{equ-main}). Conversely, if there is a radial solution to (\ref{equ-main}), using Lemma \ref{theo:ode}, we can transform it into a solution of $z''=F(t,z)$. This implies that (\ref{gk}) holds. Therefore by (\ref{klk}) we have that Condition (\ref{abc}) is true. \hfill $\square$\\


\textbf{Proof of Theorem \ref{theo:radi}(i).} Let $g(r,s)$ be as in the statement of the theorem. Because $\lim_{(|x|,s)\to(\infty,\infty)}\frac{f(x,s)}{g(x,s)}=1$, without loss of generality we can assume that $r_0$ and $s_0$ are large enough so that $g(r,s)<2f(r,s)$ on $\Omega$. Define $l(r,s):=\frac{1}{2}g(r,s)$. We work with $l(r,s)$ because we want to use the monotonicity condition (c3) which is not available for $f(|x|,u)$. We claim that if (\ref{equ-main}) has a solution, then the problem
\begin{align} \label{equ:l}
\begin{cases}
\Delta y(x) + \nabla h(x) \cdot \nabla y(x)= l(|x|,u(x)) ~~~~~~&  x \in \mathbb{R}^n ~~~ (n\geq 2),\\
y(x)\longrightarrow \infty  & x\rightarrow \infty .
\end{cases}
\end{align}
has a radial solution. Assume for the contrary that (\ref{equ-main}) has a solution while there is no radial solution to (\ref{equ:l}). To reach to a contradiction, we study another related PDE. Consider constants $s_1 > s_0$ and $r_1 > r_0$ such that
\begin{equation}\label{hhh}
\max_{|x|=r_0} u(x)< s_1 <\max_{|x|=r_1} u(x).
\end{equation}
These constants exist because $\lim_{|x|\to \infty}u(x)=\infty$. We want to  prove that $\exists r_2>r_1$ such that there exists a radial solution to the following PDE
\begin{equation}\label{equ:v}
\begin{aligned}
	\begin{cases}
		\Delta v(x) + \nabla h(x) \cdot \nabla v(x) = l(x,v)~~~& r_0<|x|<r_2\\
		v(x)=s_1 & |x|=r_0\\
		v(x)=s_1 & |x|=r_1\\
		v(x) \longrightarrow \infty & |x| \longrightarrow r_2^-.
	\end{cases}
\end{aligned}  
\end{equation}
Setting $t=p(r)$, $z(t)=v(p^{-1}(t))$, and $F(t,z)=(p^{-1}(t))^{2n-2}e^{2a(p^{-1}(t))}l(p^{-1}(t),z)$ the problem of finding $r_2$ is equivalent to finding $t_2 \in (t_1,0)$ such that there exists a solution to
\begin{equation}\label{kkk}
	\begin{cases}
		z''(t) = F(t,z)\\
		z(t_0)=z(t_1)=s_1\\
		\lim_{t \to t_2^-} z(t)=\infty .
	\end{cases}  
\end{equation}
Note that because we assumed (\ref{equ:l}) has no solution, proof of part (iii) implies that 
\begin{equation} \label{luk}
\int_ \Omega -p(|x|)e^{h(|x|)}l(x,s)dx = \infty,
\end{equation}
which again part (iii) results in
\begin{equation}\label{ccc}
-\int_{t_0}^0 tF(t,s)dt=\infty.
\end{equation}
On the bounded interval $[t_0,t_1]$ with bounded boundary values $z(t_0)=z(t_1)=s_1$, we can use the Green's function of $\frac{-d^2}{dt^2}$ to find a solution to $z''=F(t,z)$ on this domain. Let $t_2 > t_1$ be the maximal time where $z(t)$ is continuously solves $z''=F(t,z)$. Since $z''(t)=F(t,z)\geqslant 0$ and $z(t_0)=z(t_1)$, we have that $z'(t)\geqslant 0$. There are only three possibilities. The first case is when $t_2=0$ and $\lim_{t \to 0^-}z(t)=\infty$. This possibility is ruled out because (\ref{ccc}) implies that (\ref{equ:odee}) has no solution. The second possibility is that $t_2=0$ and $\lim_{t \to 0^-}z(t)<\infty$. In this case by integrating $z''=F(t,z)$ twice we have $-\int_{t_1}^0 t F(t,z)dt=z(0^-)-z(t_1)+t_1z'(t_1)< \infty$ which is a contradiction by (\ref{ccc}). The only possibility is that $t_2 \in (t_1,0)$ and $\lim_{t \to {t_2}^-}z(t)=\infty$.
Therefore we found $t_2$ with the required conditions. By converting (\ref{kkk}) back into the corresponding PDE, there exists $r_2=p^{-1}(t_2) \in (r_1,\infty)$ such that there is a radial solution to (\ref{equ:v}). The set $\Sigma =\lbrace x \in (r_0,r_2) | u(x) > v(x)  \rbrace$ is an open and non-empty because of the definition of $s_1$. Since $f(r,s)\geqslant h(r,s)$ on $\Sigma \subset \Omega$, we have
\begin{equation*}
\Delta(u-v)(x)+\nabla h(x) \cdot \nabla (u-v)(x)=f(x,u)-l(x,v)>0 ~~ \forall x \in \Sigma . 
\end{equation*}
But $u(x)-v(x)=0$ on $\partial \Sigma$. This is a contradiction by the maximum principle. Hence assumption (\ref{luk}) is not true. Therefore if Problem (\ref{equ-main}) has a solution, then $\int_ \Omega -p(|x|)e^{h(|x|)}l(x,s)dx < \infty $. Because $\lim_{(r,s) \to (\infty,\infty)}\frac{f(r,s)}{l(r,s)}=1/2$ we have
\begin{equation*}
\exists \hat{s}\geqslant s_0 ~~ \int_{|x|>r_0} -p(|x|)e^{h(|x|)}f(x,s)dx < \infty ~~~~ \forall s> \hat{s}.
\end{equation*} \hfill $\square$\\


\textbf{Proof of Theorem \ref{theo:radi}(ii).} We start by showing that that the difference of any two $C^2$ solutions $u_a(x)$ and $u_b(x)$ of PDE (\ref{equ-main}) goes to zero at infinity. First assume that $y_a(x)$ and $y_b(x)$ are two radial solutions of the PDE. By setting
\[t=p(r)~,~z(t)=y(p^{-1}(t))~,~F(t,z):=p^{-1}(t)^{2n-2}e^{2h(p^{-1}(t))}f(p^{-1}(t),z),\]
Lemma \ref{theo:ode} implies that we can find two solutions $z_a(p(|x|))=y_a(x)$ and $z_b(p(|x|))=y_b(x)$ of the corresponding ODE. By Lemma \ref{lem:integ} the difference of any two large solutions of the ODE $z''=F(t,z)$ goes to zero as $t \to 0^-$. This implies that $\lim_{|x| \to \infty}|y_a(x)-y_b(x)|=0$.\\

Let $r>0$ be large enough so that $u(x) > s_0$ for $|x|>r$. Let $m=\min_{|x|=r} u_a(x)$ and $M=\max_{|x|=r} u_a(x)$. Now consider the sequences $u_k(x)$ and $U_k(x)$ described in the remark of Lemma \ref{theo:ode}. By the construction $u_a(x)-u_k(x)>0$ on $|x|=r$ and $\lim_{|x| \to \infty} u_a(x) - u_k(x)= \infty$. Since $f_s(r,s)\geqslant 0$ we have
\begin{equation*}
\Delta (u_a(x)- u_k(x)) + \nabla h(x) \cdot \nabla (u_a(x)-u_k(x)) = f(x,u_a(x))-f(x,u_k(x))\geqslant 0 .
\end{equation*}
Therefore maximum principle implies $u_a(x)\geqslant u_k(x)$ for all $k$ and all $|x|>r$. Hence 
\[u_a(x)\geqslant u_0(x)=\lim_{k \to \infty} u_k(x)~~~\text{for}~|x|>r.\] 
Similarly $u_a(x)\leqslant U_k(x)$ for $r<|x|<p_k$. Since $\lim_{k \to \infty}p_k=\infty$ we have $u_a(x) \leqslant U_0(x)$ on $|x|>r$. Furthermore $u_0(x)$ and $U_0(x)$ are two radial solutions of the problem (\ref{equ-main}). By the discussion at the beginning of this step $\lim_{|x| \to \infty} |U_0(x)-u_0(x)|=0$. Since $u_0(x) \leqslant u_a(x) \leqslant U_0(x)$ we have $\lim_{|x| \to \infty} |u_a(x)-u_0(x)|=0$. By the similar argument for $u_b(x)$ we have $\lim_{|x| \to \infty} |u_b(x)-u_0(x)|=0$.\\

Now assume that $R$ is an orthonormal transformation on $\mathbb{R}^n$. We have
\begin{align*}
[\nabla (h(Rx))] \cdot [\nabla (u(Rx))]&=[\nabla (h(Rx))]^T [\nabla (u(Rx))]\\
&= [(\nabla h)(Rx))] R R^T [(\nabla u) (Rx)]\\
&=\nabla h \cdot \nabla u (R(x))
\end{align*}
Furthermore the Laplace operator is interchangeable with orthonormal operators in the sense that for $u_R(x)=u(R(x))$ we have $\Delta u(R(x))= \Delta u_R(x)$. Therefore for a given solution $u(x)$ of Problem (\ref{equ-main}) we have
\begin{equation*}
\Delta (u_R-u)(x)+ \nabla h \cdot \nabla (u_R-u)(x)=f(x,u_R)-f(x,u).
\end{equation*}
By the argument at the beginning of the proof we know that $\lim_{|x| \to \infty} |u_R-u|=0$. Because $f_s(r,s)\geqslant 0$, the maximum principle implies that $u_R\equiv u$. Therefore $u(x)$ is radial. \hfill $\square$\\


\section*{Appendix}

In this appendix we gather the statements of ODE lemmas required for our arguments. See \cite{T1} for the proofs of the lemmas.

\begin{lemma}\label{lem:integ}
Let $\Gamma = \lbrace (t,z)|\hat{t} \leqslant t<0,~ 0 <\hat{z} <z \rbrace$ be given. Assume that $F(t,z)$ and $F_z(t,z)$ are $C^0$ and non-negative on $\Gamma$. Assume also that $F_z(t,z)$ is superlinear in $z$ and $F(t,z)$ is monotone in $t$ on $\Gamma$. Then the problem 
\begin{equation}\label{a}
  \begin{cases}
   z''(t)=F(t,z(t)) \\
   \lim_{t \to 0^-} z(t)=\infty
  \end{cases}
\end{equation}
has a $C^2$ solution if and only if there exists $c \in (\hat{t},0)$ such that
\begin{equation}\label{aa}
-\int_{c}^0 tF(t,z)dt < \infty ~~~~ \forall z>\hat{z}.
\end{equation}
Furthermore, for any pair of solutions $z_1(t)$, $z_2(t)$ to (\ref{a}) we have 
\[\lim_{t \to 0^-} |z_2(t)-z_1(t)|=0.\]
\end{lemma}


\begin{lemma}\label{lem:v}
Let $\Gamma = \lbrace (t,z)|\hat{t} \leqslant t<0,~ 0 <\hat{z} <z \rbrace$ be given. Assume that $F(t,z)$ and $F_z(t,z)$ are $C^0$ and non-negative on $\Gamma$. Assume also that $F(t,z)$ is superlinear in $z$ on $\Gamma$. Then for each $\bar{z}>\hat{z}$ and $\bar{t}>\hat{t}$ there exists a sequence $\lbrace \rho_k \rbrace^{\infty}_{k=0} \subseteq (\bar{t},0)$ with $\lim_{k \to \infty}\rho_k = 0$ and two sequence of $C^2$ functions $\lbrace z(t) \rbrace^{\infty}_{k=o}$ and $\lbrace Z(t) \rbrace^{\infty}_{k=o}$ such that
\begin{itemize}
\item[(i)] $Z_0(t)$ and $z_0(t)$, $z_1(t)$,... are solutions of  
\begin{equation*}
  \begin{cases}
   z''(t)=F(t,z(t)) ~~~ t\geqslant t_1 \\
   z(\bar{t})=\bar{z}.
  \end{cases}
\end{equation*}
\item[(ii)]$\forall k \geqslant 1$, $Z_k(t)$ is a solution of  
\begin{equation*}
  \begin{cases}
   Z_k''(t)=F(t,Z_k(t)) ~~~ t_1 \leqslant t < \rho_k \\
   \lim_{t \to \rho_k^-} Z_k(t)=\infty \\
   Z_k(\bar{t})=\bar{z}.
  \end{cases}
\end{equation*}
\item[(iii)]$\lim_{t \to 0^-} z_0(t) =\lim_{t \to 0^-} Z_0(t)=\infty$.
\item[(iv)]$\forall k \geqslant 1$, $z_k(t)$ is finite as $t \to 0^-$.
\item[(v)] For each $t \in (\bar{t},0)$ we have 
\[\lim_{k \to \infty}z_k(t)=z_0(t) ~~~ and ~~~ \lim_{k \to \infty}Z_k(t)=Z_0(t).\]
\end{itemize}
\end{lemma}


\bibliographystyle{plain}
\bibliography{radial}

\end{document}